# Asymptotic spectral flow for Dirac operators

C. H. Taubes[†]

Department of Mathematics
Harvard University
Cambridge, MA  02138

chtaubes@math.harvard.edu

Let M denote a compact, oriented Riemannian manifold of odd dimension $n \geq 3$. Suppose that $F \to M$ is a principle bundle with structure group $Spin(n) \times_{\{\pm 1\}} U(k)$ such that $F/U(k)$ is the priniciple $SO(n)$ bundle of orthonormal frames for TM. A connection on the principle bundle $F/Spin(n) \to M$ determines a self-adjoint Dirac operator on a certain associated Clifford module. This understood, suppose that $\{A_s\}_{s \in [0,1]}$ is a differentiable path of such connections. An estimate is given here for the spectral flow along the corresponding 1-parameter family of Dirac operators.

[†]Supported in part by the National Science Foundation

# 1. Introduction

This note generalizes certain estimates from [T] for the spectral flow for a 1-parameter family of Dirac operators. In this regard, the estimates that appear in [T] concern only $Spin_{\mathbb{C}}$ Dirac operators on 3-dimensional manifolds, whereas the techniques that prove these estimates provide analogous estimates for the spectral flow families of $Spin_{\mathbb{C}}$ Dirac operators on any odd dimensional compact manifold, and also for twisted versions of such operators. These generalizations were not discussed in [T] as the latter addressed a particular application. This note is meant to highlight the generalizations.

To say more, suppose that $n \geq 1$ is odd and that M is a compact, oriented n-dimensional Riemannian manifold. Let $Fr \to M$ denote the bundle of oriented, orthonormal frames. Let $k \geq 0$ and suppose that this bundle lifts to a $Spin(n) \times_{\mathbb{Z}/2} U(k)$ bundle, $F \to M$. Here $\mathbb{Z}/2$ acts on $Spin(n)$ as its center, and it acts on $U(k) \subset GL(k; \mathbb{C})$ as + and - times the identity matrix. Let S denote the fundamental spinor representation of $Spin(n)$, this a complex vector space of dimension $2^{(n-1)/2}$. With S understood, introduce the complex hermitian vector bundle $\mathbb{S} \to M$ that is associated to F via the representation $S \otimes_{\mathbb{C}} \mathbb{C}^k$ of $Spin(n) \times_{\mathbb{Z}/2} U(k)$.

The Levi-Civita connection on Fr and a unitary connection on the principle $PU(k) = U(k)/\{\pm 1\}$ bundle $P = F/Spin(n)$ together give $\mathbb{S}$ a unitary connection. With a unitary connection, A, on P chosen, use $D_A: C^{\infty}(M; \mathbb{S}) \to C^{\infty}(M; \mathbb{S})$ to denote the resulting Dirac operator. The latter is obtained by composing the covariant derivative, $\nabla_A: C^{\infty}(M; \mathbb{S}) \to C^{\infty}(M; \mathbb{S} \otimes T^*M)$ with the endomorphism from $\mathbb{S} \otimes T^*M$ to $\mathbb{S}$ which is dual to the Clifford multiplication endomorphism from $T^*M$ to $End(\mathbb{S})$. This Dirac operator extends to $L^2(M; \mathbb{S})$ as an unbounded, self-adjoint operator with domain $L^2_1(M; \mathbb{S})$. It has purely point spectrum and each eigenvalue has finite multiplicity. Moreover, the set of eigenvalues is unbounded from below and above, and has no accumulation points.

Let $A_0$ and $A_1$ denote a pair of unitary connections on P such that neither of the $A = A_0$ and $A = A_1$ versions $D_A$ have eigenvalue 0. The Dirac spectral flow from $A_0$ to $A_1$ is defined as follows: Fix a path of unitary connections on P that starts at $A_0$ and ends at $A_1$. The Dirac spectral flow constitutes an algebraic count of the connections on this path where the corresponding Dirac operator has a zero eigenvalue. For those unfamiliar with the notion of spectral flow, here is a brief definition (more is said below): Note first that the eigenvalues vary continuously along the path. This understood, a zero crossing of a non-degenerate eigenvalue contributes to the count with +1 weight if the eigenvalue crosses zero from a negative to a positive value as s increases. On the other hand, a zero crossing of a non-degenerate eigenvalue contributes to the count with -1 if the eigenvalue crosses zero from a positive to a negative value as s increases. If the path is suitably generic, then only these two cases arise. The resulting count of the crossings turns out to be path independent and so depends only on the ordered pair $(A_0, A_1)$. This count is the Dirac spectral flow function; it is denoted in what follows by $f(A_0, A_1)$.

The statement below of the main theorem refers to a differential form on M that is constructed from the Riemannian curvature tensor. This form is denoted by $\Omega_{\hat{A}}$; it is the canonical closed form that is constructed using the Riemannian curvature tensor so as to represent the $\hat{A}$-genus of M in deRham cohomology. To elaborate, let $\mathfrak{so}(n)$ denote the Lie algebra of SO(n). Following the exposition in Chapter 1 of [BGV], define a map from $\mathfrak{so}(n)$ to itself by using the convergent power series for $\sigma \to \frac{\sinh(\sigma)}{\sigma}$. Now, let $\mathcal{R}$ denote the Riemann curvature but viewed as a 2-form on M with values in $\text{Fr} \times_{ad} \mathfrak{so}(n)$. Then $\mathfrak{j} = \det(\frac{\sinh(\mathcal{R}/2)}{\mathcal{R}/2})$ defines a differential form on M with 0-form component equal to 1. This understood, then $\mathfrak{j}^{-1/2}$ also defines a differential form on M with 0-form component equal to 1. The form $\Omega_{\hat{A}}$ is defined to be $\mathfrak{j}^{-1/2}$.

The main theorem also refers to a relative Chern-Simons form; in the present context, this is a certain differential form on M that is determined by an ordered pair of connections on P. To set the stage for the definition, suppose that A is a connection on P. Let $A_F$ denote the pull-back of A to F, this a 1-form on F with values in the Lie algebra of U(k). Use $F_{A_F}$ to denote the curvature 2-form of $A_F$. This 2-form is obtained from the curvature 2-form of A by dividing the central component of the latter by 2. With this notation set, let $(A_0, A_1)$ denote a given ordered pair of connections on P. Write $A_{1F} = A_{0F} + \hat{a}$, where $\hat{a}$ is a 1-form on M with values in the bundle that is associated to P by the action of U(k) on the Lie algebra of U(k). For $s \in [0, 1]$, introduce $A_F(s) = A_{0F} + s\hat{a}$. The relative Chern-Simons form for the ordered pair $(A_0, A_1)$ is

$$\mathfrak{chs}(A_0, A_1) = \int_0^1 \text{tr}_{\mathbb{C}^k}(\hat{a} \wedge \exp(F_{A_F(s)})) \, ds \,.$$

(1.1)

This form has the following useful property: Let $\mu$ denote a closed form on P. Then $\int_M (\mu \wedge \mathfrak{chs}(A_0, A_1))$ is not changed if the path of connections $s \to A_{F(s)}$ and the 1-form $\hat{a}$ (which is $\frac{d}{ds} A_F(s)$) that appear in (1.1) are simultaneously replaced by $A_{sF}$ and $\frac{d}{ds} A_{sF}$ where $s \to A_s$ is any piece-wise differentiable path of connections that begins at $A_0$ when $s = 0$ and ends at $A_1$ when $s = 1$.

One more bit of notation and one additional remark are needed to set the stage for the upcoming theorem. Here is the extra notation: If A is a connection on P, the theorem uses r(A) to denote the infimum of numbers $r \in [1, \infty)$ such that

$$r^{-1} |F_A| + r^{-3/2} |\nabla_A F_A| + \cdots + r^{-n/2} |(\nabla_A)^{(n-1)/2} F_A| \leq 1 \,.$$

(1.2)

The promised remark concerns a certain subgroup of the group of automorphisms of the bundle P. To set the stage, note that the group U(k) acts on itself via conjugation and this action factors to an action of PU(k) on U(k). The latter action is denoted in what follows by Ad. Any section of the bundle $P \times_{Ad} U(k)$ determines an automorphism of P, and all autormorphisms of this type lift to give automorphisms of F that cover the identity automorphism of the frame bundle Fr. When g denotes an automorphism of P and A denotes a connection on P, then g·A is used to denote the pull-back of A by the automorphism.

What follows is the principle result of this note:

**Theorem 1:** *There exist constants, $\kappa \geq 1$ and $p \in (\frac{n-3}{2}, \frac{n-1}{2})$ with the following significance: Let $\{A_s\}_{s \in [0,1]}$ denote a differentiable path of connections on P and set r to denote the maximum from the set $\{r(A_s)\}_{s \in [0,1]}$. Let g denote a section of $P \times_{Ad} U(k)$. Then $f(A_0, g \cdot A_1)$ differs from*

$$(\tfrac{1}{2\pi i})^{(n+1)/2} \int_M (\Omega_{\hat{A}} \wedge \mathfrak{chs}(A_0, g \cdot A_1))$$

(1.3)

*by at most $\kappa(1 + r^{p+1} + \kappa \int_0^1 (r(A_s)^p (\int_M |\tfrac{d}{ds} A_s|)) \, ds)$.*

Neither the expression in (1.3) nor the spectral flow function are gauge invariant except under certain particular situations. To elaborate, suppose that g is a section of the bundle $P \times_{Ad} U(k)$. Then the spectrum of $D_{g \cdot A}$ and $D_A$ are identical. Even so, the Dirac spectral flow from A to g·A need not vanish. However, it follows from what is done in the seminal paper by Atiyah, Patodi and Singer [APS] that $f(A, g \cdot A)$ is identical to the $A_0 = A$ and $A_1 = g \cdot A$ version of (1.3). Thus, the difference between $f$ and the the expression in (1.3) is, in fact, gauge invariant.

Although the expression in (1.3) and the spectral flow depend only on the endpoints of the path $s \to A_s$, the error term given by the theorem depends on the chosen path. It is not clear whether there exists a useful apriori bound on the error term for *some* path between any two given connections $A_0$ and $A_1$ in the case of U(k) with k > 1 that depends only on $A_0$ and $A_1$. Results from [U] may have some bearing on this issue.

Here is what can be said when k = 1:

**Theorem 2**: *Suppose that F is a $\mathrm{Spin}_{\mathbb{C}} = \mathrm{Spin}(n) \times_{\{\pm 1\}} U(1)$ lift of the frame bundle of M. There exists $\kappa \geq 1$ and $p \in (\frac{n-1}{2}, \frac{n+1}{2})$ with the following significance: Let $A_0$ and $A_1$ be a pair of connections on P = F/Spin(n). Set $r = \max\{1, r(A_0), r(A_1)\}$. Then $f(A_0, A_1)$ differs from $(\tfrac{1}{2\pi i})^{(n+1)/2} \int_M (\Omega_{\hat{A}} \wedge \mathfrak{chs}(A_0, A_1))$ by at most $\kappa \, r^p$.*

Theorem 2 suggests a question: Let P be as in Theorem 2, let $A_0$ denote a connection on P and let a denote an imaginary number valued 1-form on M. What is the subleading order behavior of the spectral flow function $f(A_0, A_0 + ra)$ as $r \to \infty$? Note that Theorem 2 asserts that the leading order term is

$$r^{(n+1)/2} \left(\frac{1}{4\pi i}\right)^{(n+1)/2} \frac{1}{(\frac{1}{2}(n+1))!} \int_M (a \wedge (da)^{(n-1)/2})$$

(1.4)

Theorems 1 and 2 have a corollary that follow using the results [APS]. As background for the corollary, introduce the manifold $X = \mathbb{R} \times M$. In what follows, s is used to denote the standard Euclidean coordinate on $\mathbb{R}$. Give X the product Riemannian metric. The projection from X to M pulls back $\mathbb{S}$ to give a vector bundle over X, this also denoted by $\mathbb{S}$.

Fix a pair of connections, $A_0$ and $A_1$, on P, and fix a differentiable path $s \to \mathcal{A}(s)$ from $\mathbb{R}$ to the space of connections on P such that $\mathcal{A}(s) = A_0$ for $s \ll -1$ and $\mathcal{A}(s) = A_1$ for $s \gg 1$. With $\mathcal{A}$ chosen, define the operator $\mathcal{D}: C^\infty(X; \mathbb{S}) \to C^\infty(X; \mathbb{S})$ by

$$\mathcal{D} = \frac{\partial}{\partial s} + D_{\mathcal{A}(\cdot)} .$$

(1.5)

The operator $\mathcal{D}$ defines a Fredholm operator from $L^2_1(X; \mathbb{S})$ to $L^2(X; \mathbb{S})$. As such, $\mathcal{D}$ has an index, and this index is given, according to (4.3) in [APS] by the following formula:

$$\text{index}(\mathcal{D}) = \left(\frac{1}{2\pi i}\right)^{(n+1)/2} \int_M (\Omega_{\hat{\mathcal{A}}} \wedge \mathfrak{chs}(A_0, A_1)) - \tfrac{1}{2}(\eta_{A_1} - \eta_{A_0})$$

(1.6)

where the notation is such that $\eta_A$ denotes the spectral asymmetry function of $D_A$ when A is a connection on P. This spectral asymmetry function is defined as follows: It is the value at 0 of the analytic continuation to $\mathbb{C}$ of the function on the $s \gg 1$ part of the real line that sends s to

$$\sum_\varsigma \text{sign}(\lambda_\varsigma) |\lambda_\varsigma|^{-s} ;$$

(1.7)

here, the sum is indexed by an orthonormal basis of eigenfunctions of $D_A$ and $\lambda_{(\cdot)}$ denotes the eigenvalue of the indicated eigenvector. Theorem 3.10 in [APS] asserts that this analytic continuation is finite at $s = 0$.

As noted in [APS], the index of $\mathcal{D}$ is precisely $f(A_0, A_1)$. Thus, (1.6) and Theorems 1 and 2 have the following to say about $\eta_{A_1} - \eta_{A_0}$:

**Corollary 3:** *There exist constants $\kappa \geq 1$ and $p \in (\frac{n-3}{2}, \frac{n-1}{2})$ with the following significance: Let $\{A_s\}_{s \in [0,1]}$ denote a differentiable path of connections on P. Set r to denote the maximum from the set $\{r(A_s)\}_{s \in [0,1]}$. Then the absolute value of the difference between the respective spectral asymmetry functions of $A_0$ and $A_1$ obeys*

$$|\eta_{A_1} - \eta_{A_0}| \leq \kappa (1 + r^{p+1} + \int_0^1 (r(A_s))^p (\int_M |\tfrac{d}{ds} A_s|)) \, ds ) \, .$$

*In the case where $k = 1$, there exists $p' \in (\frac{n-1}{2}, \frac{n+1}{2})$ such that the following is true: Let $A_0$ and $A_1$ denote a pair of connections on P, and set r to equal the maximum of 1, $r(A_0)$ and $r(A_A)$. Then $|\eta_{A_1} - \eta_{A_0}| \leq \kappa \, r^{p'}$.*

The strategy for the proofs of Theorems 1 and 2 is the same as that used in Section 5 of [T] to prove a 3-dimensional version of Theorem 2. This strategy is outlined next, and the analysis follows.

### a) Spectral flow

This and the next subsection review some constructions from Sections 5a and 5b in [T]. To start, suppose that $\{A_s\}_{s \in [0,1]}$ is a family of connections on P with a real analytic parametrization by [0, 1]. The spectral flow for the family $\{D_{A_s}\}_{s \in [0,1]}$ is defined with the help of a certain stratified, real-analytic set in $\mathbb{R} \times [0, 1]$. This set is denoted by $\mathcal{E}$, and its stratification is depicted by

$$\mathcal{E} = \mathcal{E}_1 \supset \mathcal{E}_2 \supset \cdots .$$

(1.8)

where $\mathcal{E}_k$ consists of the set of pairs $(\lambda, s)$ such that $\lambda$ is an eigenvalue of $D_{A_s}$ with multiplicity k or greater. Each $\mathcal{E}_k$ is a closed set. Moreover, as can be proved using results in Chapter 7 of [K], each $\mathcal{E}_{k*} = \mathcal{E}_k - \mathcal{E}_{k+1}$ is an open, real analytic submanifold in $\mathbb{R} \times [0, 1]$. As in [T], the collection $\{\mathcal{E}_{k*}\}$ are called the *smooth strata* of $\mathcal{E}$. When the 1-dimensional smooth strata are oriented by the pull-back from $\mathbb{R} \times [0, 1]$ of the 1-form ds, then the zero dimensional strata can be consistently oriented so that the formal, weighted sum $\mathcal{E}_* = \mathcal{E}_{1*} + 2\mathcal{E}_{2*} + \cdots$ defines a locally closed cycle in $\mathbb{R} \times [0, 1]$. This also follows from what is said in Chapter 7 of [K]. This means the following: Let f denote a smooth function on $\mathbb{R} \times (0, 1)$ with compact support. Then

$$\Sigma_{k=1,2,\ldots} k \int_{\mathcal{E}_{k*}} df = 0$$

(1.9)

Sard's theorem finds a dense, open set $\mathbb{U} \subset \mathbb{R}$ with the property that the respective maps from a point, $*$, to $\mathbb{R} \times [0, 1]$ that send $*$ to $(\lambda, 0)$ and to $(\lambda, 1)$ are transversal to the smooth strata of $\mathcal{E}$ for all $\lambda \in \mathbb{U}$. In this language, the spectral flow for the family $\{\mathcal{L}_s\}_{s \in [0,1]}$ is defined as follows: Fix $\lambda_0 \in \mathbb{U}$ with $\lambda_0 > 0$. By Sard's theorem, there exist smooth, oriented paths $\sigma \subset \mathbb{R} \times [0, 1]$ that start at $(\lambda_0, 0)$, end at $(\lambda_0, 1)$, and are transversal to the smooth strata of $\mathcal{E}$. Since both the $A = A_0$ and $A = A_1$ versions of $D_A$ do not have zero as an eigenvalue, it follows that $0 \in \mathbb{U}$. This understood, choose a path $\sigma$ that starts at $(0, 0)$, ends at $(0, 1)$ and is transversal to the smooth strata of $\mathcal{E}$. Such a path has a well defined intersection number with $\mathcal{E}$, this being

$$f = \sum_{k=1,2,\ldots} \sum_{p \in \sigma \cap \mathcal{E}_{k*}} (-1)^{o(p)} k ,$$

(1.10)

where $o(p) \in \{0, 1\}$. In the case where $\sigma$ is the graph of a smooth function from $[0, 1]$ to $\mathbb{R}$, the sign $(-1)^{o(p)}$ is obtained as follows: The pull-back to a smooth, 1-dimensional stratum of $\mathcal{E}$ of the 1-form $d\lambda$ from $\mathbb{R} \times [0, 1]$ at a point $(\lambda, s)$ can be written as $\lambda' ds$ with

$$\lambda' = \langle \varsigma, \mathrm{cl}(\tfrac{d}{ds} A_{sF}) \varsigma \rangle_{L^2} .$$

(1.11)

Here, the notation uses $\varsigma$ to denote an eigenvector of $D_{A_s}$ with $L^2$ norm 1 whose eigenvalue is $\lambda$, and $\langle\,,\,\rangle_{L^2}$ denotes the inner product on $L^2(M; \mathbb{S})$. The sign of $\lambda'$ at an intersection point with the image of a graph is the factor $(-1)^{o(\cdot)}$ that appears in (1.11).

**b) Estimating spectral flow**

With $f$ understood, what follows explains how Theorem 1's estimate for $f$ is derived. To this end, introduce the constant r from Theorem 1. Fix $t \in (0, r^{-1})$; a specific choice is made in Section 3. With t chosen, define the orientation preserving diffeomorphism $\Phi$ from $\mathbb{R}$ to $(-(\tfrac{\pi}{t})^{1/2}, (\tfrac{\pi}{t})^{1/2})$ by the formula

$$\Phi(\lambda) = \int_0^\lambda e^{-\rho^2 t} d\rho .$$

(1.12)

Fix $R \geq 1$ and set $T = \Phi(Rt^{-1/2})$. A specific choice for R is also made Section 3. Note for reference later that

$$|t^{1/2} T - (\tfrac{\pi}{4})^{1/2}| \leq \tfrac{1}{2R} e^{-R^2} .$$

(1.13)

Now let S denote the circle that is obtained from the interval [-T, T] by identifying the endpoints. This circle has a fiducial point, $T_*$, that given by $\{\pm T\}$, and an orientation given by the orientation of (-T, T).

For each $s \in [0, 1]$, let $\mathfrak{n}_s$ denote the maximal number of linearly independent eigenvectors of $D_{A_s}$ whose eigenvalue lies in $[-Rt^{-1/2}, Rt^{-1/2}]$. Use $\mathfrak{n}$ in what follows to denote the maximum from the set $\{\mathfrak{n}_s\}_{s \in [0,1]}$. An estimate for the spectral flow $f$ is obtained by considering the trajectories of $\mathfrak{n}$ particles on S whose paths vary continuously and piecewise differentiably as functions of $s \in [0, 1]$.

To elaborate, introduce $\mathcal{E}^R$ to denote the set $\{(\lambda, s) \in \mathcal{E}: |\lambda| < Rt^{-1/2}\}$, and for each k, use $\mathcal{E}^R_{k*}$ to denote $\mathcal{E}_{k*} \cap \mathcal{E}^R$. Each point $(\lambda, s) \in \mathcal{E}^R_{k*}$ corresponds to k particles on S all at the point $\Phi(\lambda)$. If $\mathcal{E}^R_{k*}$ is 1-dimensional, then these k particles all move together near s, and the common tangent vector to their trajectories is $\lambda'(\frac{d}{d\lambda}\Phi)|_\lambda$ with $\lambda'$ as in (1.11). The set of all such trajectories that limit to a given zero-dimensional stratum, $\mathcal{E}^R_{k'*}$ as s limits to some $s_*$ can be joined at this stratum to obtain a set of k′ continuous, piecewise smooth, oriented trajectories that are defined for s near $s_*$. This follows from (1.9). There is no canonical way to do this joining, but any method will suffice.

At any given value of s, what was just described accounts for at most $\mathfrak{n}_s$ of the particles. The remaining particles are at the point $T_* \in S$. Particles move off or onto the point $T_*$ at values of s for which either of the points $(-Rt^{-1/2}, s)$ or $(Rt^{1/2}, s)$ are in the closure of $\mathcal{E}^R$. The particles that move on or off $T_*$ and the direction in S that they move are determined by which smooth strata of $\mathcal{E}^R$ have $(-Rt^{-1/2}, s)$ or $(Rt^{1/2}, s)$ in their closure. The rules for this are essentially identical to those given in the preceding paragraph.

Granted the preceding, let $s \to z(s) \in S$ denote the trajectory of a given particle. Let $\Delta z$ denote $z(1) - z(0)$; this is the net change in z as s increases from 0 to 1. The intersection number with the point $0 \in S$ of this trajectory is, at most, the least integer that is greater than $\frac{1}{2T}\Delta z$, thus at most $\frac{1}{2T}\Delta z + 1$. Meanwhile, this intersection number is at least the greatest integer less than $\frac{1}{2T}\Delta z$, thus at least $\frac{1}{2T}\Delta z - 1$. In this regard, note that where the trajectory $s \to z(s)$ is differentiable, the chain rule finds its derivative to be

$$\tfrac{d}{ds} z = \lambda'\,(\tfrac{d}{d\lambda}\Phi)(\lambda) = \lambda' e^{-\lambda^2 t}$$

(1.14)

where $\lambda = \Phi^{-1}(z)$ is the corresponding eigenvalue and where $\lambda'$ is given by (1.11).

As a consequence of what was just said, the spectral flow $f$ differs by at most $\mathfrak{n}$ from the integral between 0 and 1 of the function

$$\wp(s) = \tfrac{1}{2T}\sum_{\varsigma \in \Theta}\left(\int_M \varsigma^\dagger \mathrm{cl}(\tfrac{d}{ds} A_{sF})\varsigma\right) e^{-\lambda_\varsigma^2 t},$$

(1.15)

where the sum is over a basis, $\Theta$, of orthonormal eigenvectors of $D_{A_s}$ whose eigenvalue has absolute value no greater than $Rt^{-1/2}$, and where the notation again uses $\lambda_{(\cdot)}$ to denote the eigenvalue of the indicated eigenvector. With regards to the integral of the function that sends s to $\wp(s)$, note that $\wp$ varies with s as a piecewise continuous and bounded function.

The estimates in Theorem 1 for the spectral flow are obtained by deriving suitable estimates for the function $s \to \wp(s)$ and upper bounds for the number $\mathfrak{n}$. Both of these tasks are accomplished with the help of the heat kernel for the operator $(D_{A_s})^2$. The next section supplies what is required from the heat kernel. The results from Section 2 are put to use in Section 3 to give proofs of Theorems 1 and 2.

## 2. The heat kernel

The heat kernel for the square of the Dirac operator proves its worth in two ways. First, it provides an upper bound for $\mathfrak{n}$. Second, it provides an estimate for the function $\wp(s)$ that appears in (1.15). This subsection summarizes some basic facts about the heat kernel that are needed for these applications.

### a) Bounds on the heat kernel

To start the story, fix a connection, A, on P. Let $\pi_L$ and $\pi_R$ denote the respective projections from $(0, \infty) \times M \times M$ to the left and right hand factors of M. The heat kernel for $D_A^2$ is the smooth section $\text{Hom}(\pi_R^*\mathbb{S}, \pi_L^*\mathbb{S})$ over $(0, \infty) \times M \times M$ given by

$$E_A(t; x, y) = \sum_\varsigma \varsigma(x) \varsigma^\dagger(y) \, e^{-\lambda_\varsigma^2 t} ,$$

(2.1)

where the sum is indexed by a complete, orthonormal basis of eigenvectors for $D_A$. As a function of t and x with y fixed, the heat kernel obeys the equation

$$\tfrac{\partial}{\partial t} E_A = -D_A^2 E_A .$$

(2.2)

Furthermore, the $t \to 0$ limit of $E_t$ exists as a bundle valued measure:

$$\lim_{t \to 0} E_A(t, \cdot, y) = \mathbb{I} \, \delta_y(\cdot) ,$$

(2.3)

where $\mathbb{I}$ denotes the identity homomorphism in $\text{End}(\mathbb{S})$, and $\delta_y$ denotes the Dirac measure with mass at x. Note that integration, $\eta \to \int_M E_A(t; x, \cdot) \eta(\cdot)$ defines $E_A|_t$ as a bounded operator on $L^2(M; \mathbb{S})$ whose $t \to 0$ limit is the identity operator.

The following proposition gives estimates for the heat kernel. This proposition introduces $F_A$ to denote the curvature 2-form for A. Note that the norms of $F_A$ and $F_{A_F}$ differ by at most a factor of 2.

**Proposition 2.1:** *There exists a constant $\kappa \geq 1$ with the following significance: Let* A *denote a connection on* P *and let* $r \geq 1$ *be such that* $\sup_M |F_A| \leq r$. *Then*

$$|E_A(t; x, y)| \leq \kappa (\tfrac{1}{4\pi t})^{n/2}\, e^{-\text{dist}(x,y)^2/4t}\, e^{\kappa rt}$$

(2.4)

*Proof of Proposition 2.1*: The proof starts with the Bochner-Weitzenboch identity for the operator $D_A^2$:

$$D_A^2 = \nabla_A^\dagger \nabla_A + \tfrac{1}{4} R + cl(F_{A_F}) .$$

(2.5)

Here, R denotes the scalar curvature of M. To prove the inequality, fix y and set $h(t, x) = |E_A(t; x, y)|$. By virtue of (2.5), this function obeys the differential inequality

$$\tfrac{\partial}{\partial t} h \leq -d^\dagger dh + c_0(1 + r)h .$$

(2.6)

Here, and in what follows $c_0$ is used to denote a constant that depends only on the Riemannian metric on M and the integer k. The precise value can differ at each appearance. With (2.6) understood, set $h'(t, x) = e^{-c(1+r)t} h(t, x)$. The latter obeys

$$\tfrac{\partial}{\partial t} h' \leq -d^\dagger dh' .$$

(2.7)

It follows from (2.7) and (2.3) using the maximum principle that

$$h'(t; x) \leq c_0\, H_0(t; x, y)$$

(2.8)

where $H_0$ is the heat kernel for the Laplacian $d^\dagger d$. The bound asserted by Proposition 2.1 follows from (2.8) given the bounds from [P] or [M] for $H_0$.

**b) Upper bounds for a version of n**

The first application of the heat kernel enters into the proof of the following proposition.

**Proposition 2.2**: *There exists a constant, $\kappa \geq 1$, with the following significance: Suppose that* A *is a connection on* P *and that* $\sup_M |F_A| \leq r$. *For $\lambda > 0$, let $\mathfrak{n}_A(\lambda)$ denote the number* *of linearly independent eigenvectors of* $D_A$ *whose eigenvalue has absolute value less than* $\lambda$. *Then* $\mathfrak{n}_A(\lambda) \leq \kappa\,(\lambda + r^{1/2})^n$.

*Proof of Proposition 2.2*: The argument mimics an argument first used by Chen and Li in [CL] to obtain bounds on the number of linearly independent eigenvectors for the Laplacian $d^\dagger d$ with eigenvalue less than a given positive number. To start, note that by virtue of (2.1), the trace on $\mathrm{End}(\mathbb{S})$ of the endomorphism $E_A(t; x, x)$ obeys

$$\int_M \mathrm{tr}_{\mathbb{S}}(E_A(t;x,x)) \geq e^{-\lambda^2 t}\,\mathfrak{n}_A(\lambda)$$

(2.9)

What with Proposition 2.1, this implies that

$$\mathfrak{n}_A(\lambda) \leq c\left(\tfrac{1}{4\pi t}\right)^{n/2} e^{(\lambda^2 + cr)t}$$

(2.10)

Here, c is independent of $\lambda$, t and A. Taking $t = (\lambda^2 + c\,r)^{-1}$ gives the bound claimed by the proposition.

### c) Estimating a version of $\wp$

Let A denote a connection on P and let â denote a 1-form on M with values in the vector bundle that is associated to P by the adjoint representation of PU(k) on its Lie algebra. Fix $\lambda > 1$. This subsection explains how the heat kernel for $D_A^2$ is used to derive an estimate for

$$\mathfrak{p}(\lambda) = \sum_{\varsigma \in \Theta(\lambda)} \left( \int_M \varsigma^\dagger \mathrm{cl}(\hat{a})\varsigma \right) e^{-\lambda_\varsigma^2 t},$$

(2.11)

where $\Theta(\lambda)$ here denotes an orthonormal basis of eigenvectors of $D_A$ whose eigenvector has absolute value less than $\lambda$.

**Proposition 2.3**: *There exists a constant* $\kappa \geq 1$ *with the following significance: Let* A *denote a connection on* P *and let* $r \geq 1$ *be such that* $\sum_{0 \leq j \leq (n-1)/2} r^{-(1+j/2)}\,|(\nabla_A)^j F_A| \leq 1$. *Let* â *denote a 1-form on* M *with values in the bundle that is associated to* F *via the adjoint representation of* U(k) *on its Lie algebra. Fix* $\lambda > 1$ *and define* $\mathfrak{p}$ *as in (2.11) with* t *chosen so that* $rt \leq 1$. *Then*

$$\mathfrak{p}(\lambda) = \pi^{1/2}\left(\tfrac{1}{2\pi i}\right)^{(n+1)/2} t^{-1/2} \int_M (\Omega_{\hat{A}} \wedge \mathrm{tr}_{F \times_{\mathrm{ad}} \mathbb{C}^k}(\hat{a} \wedge \mathrm{ch}(F_{A_F}))) + \mathfrak{r},$$

*where* $|\mathfrak{r}| \leq \kappa\,(t^{1/2} r^{(n+1)/2} + t^{-n/2} e^{\kappa r t}\,e^{-\lambda^2 t}) \int_M |\hat{a}|$.

***Proof of Proposition 2.3***: The first step is to approximate $\wp$ by the expression that is obtained from the right hand side of (2.20) by removing the eigenvalue restriction for the sum. This amounts to writing

$$\mathfrak{p}(\lambda) = \int_M \text{tr}_\mathbb{S}(\text{cl}(\hat{a})|_x E_A(t;x,x)) - \sum_\varsigma (\int_M \varsigma^\dagger \text{cl}(\hat{a})\varsigma) \, e^{-\lambda_\varsigma^2 t},$$

(2.12)

where the sum on the far right is indexed by an orthonormal basis of eigenvectors of $D_A$ whose eigenvalue has absolute value no less than $\lambda$. To estimate the far right sum in (2.12), remark that it is, in any event, no greater than

$$\sup_M (\sum_\varsigma |\varsigma|^2 e^{-\lambda_\varsigma^2 t}) \int_M |\hat{a}|,$$

(2.13)

where the sum here is also indexed by an orthonormal basis of eigenvectors of $D_A$ whose eigenvalue has absolute value no less than $\lambda$. According to (2.1), the sum that appears in (2.13) is no greater than $e^{-\lambda^2 t/2} \text{tr}_\mathbb{S}(E_A(\frac{1}{2}t;x,x))$. This understood, an appeal to Proposition 2.1 bounds the sum that appears on the right hand side of (2.12) by

$$c_0 t^{-n/2} e^{\kappa rt} e^{-\lambda^2 t} \int_M |\hat{a}|.$$

(2.14)

Granted this bound for the sum that appears on the right hand side in (2.12), the next task is to obtain an estimate for the integral that appears on the right hand side of (2.12). The next proposition gives the desired estimate, and with (2.14), implies Proposition 2.3.

**Proposition 2.4**: *There exists a constant $\kappa \geq 1$ with the following significance: Let* A *denote a connection on* P *and let* $r \geq 1$ *be such that* $\sum_{0 \leq j \leq (n-1)/2} r^{-(1+j/2)} |(\nabla_A)^j F_A| \leq 1$. *Let* â *denote a 1-form on* M *with values in the bundle that is associated to* F *via the adjoint representation of* U(k) *on its Lie algebra. If* $t > 0$ *and* $rt \leq 1$, *then*

$$\text{tr}_\mathbb{S}(\text{cl}(\hat{a}|_x) E_A(t;x,x)) = \pi^{1/2} t^{-1/2} (\tfrac{1}{2\pi i})^{(n+1)/2} * (\Omega_{\hat{A}} \wedge \text{tr}_{\mathbb{C}^k}(\hat{a} \wedge \text{ch}(F_{A_F}))) + \mathfrak{r}$$

*where* $|\mathfrak{r}| \leq c_0 t^{-1/2} r^{(n-1)/2}(rt)|\hat{a}|$. *Here, * denotes the metric's Hodge dual operator.*

***Proof of Proposition 2.4***: There are well known techniques for generating small time asymptotics of the heat kernel. See, for example Chapter 4 in [BGV]. However, the standard asymptotics must be augmented with explicit bounds on the error. These bounds are derived in what follows.

To start, fix a Gaussian coordinate chart $\phi: B \to M$ and use $\phi$ to identify B with $\phi(B)$. This done, introduce the Euclidean coordinates $x = (x_1, \ldots, x_n)$ on B. Use parallel transport out from the origin along Euclidean rays via the Levi-Civita connection and via the connection A to trivialize the bundle F over B. With this trivialization understood, the connection $A_F$ appears as $A_F = \alpha_\nu \, dx_\nu$ where the convention is such that repeated indices are summed. Here, each $\alpha_\nu$ is a function on B with values in the Lie algebra of U(k), each vanishes at the origin, and $x_\nu \alpha_\nu = 0$. Note for future reference that $\{\alpha_\nu\}$ is determined by the curvature $F_A$ via

$$\alpha_\nu|_x = \int_0^1 \tau d\tau \, x^\mu (F_A)_{\mu\nu}|_{\tau x},$$

(2.15)

Granted this notation, then the operator $D_A{}^2$ can be written as

$$D_A{}^2 = \partial_\nu \partial_\nu + cl(F_{A_F}) + V,$$

(2.16)

where V has the following schematic form:

$$V = \gamma_{\nu\mu} \partial_\nu \partial_\mu + (\delta_{\nu\mu} + \gamma_{\nu\mu})(2\alpha_\nu \partial_\mu + \partial_\nu \alpha_\mu + \alpha_\nu(\alpha_\mu + \gamma_\mu) + 2\Gamma_\nu \partial_\mu) + \Gamma_0.$$

(2.17)

Here, $\{\gamma_{\nu\mu}, \gamma_\nu\}$ are functions on B with $|\gamma_{\nu\mu}| \le c_0 |x|^2$ and $|\gamma_\nu| \le c_0 |x|$. To say more about $\{\Gamma_\nu\}$ and $\Gamma_0$, remark that they are functions on B with values in the endomorphisms of the spinor bundle $\mathbb{S}$. Furthermore, both lie in the image of the extended Clifford homomorphism on $\Lambda^* T^* M$ of the even forms. In particular, $\Gamma_0$ lies in the image of the forms of degree 4 or less, and each $\Gamma_\nu$ lies in the image of the forms of degree 2 or less. They depend only on the Riemannian metric, and are such that $|\Gamma_\nu| \le c_0 |x|$ and $|\Gamma_0| \le c_0$.

Fix a smooth function, $\chi$, from $[0, \infty)$ to $[0, 1]$ that has value 1 on $[0, \frac{1}{16}]$ and vanishes on $[\frac{1}{8}, \infty)$. Let $\chi_\rho$ denote the function on $\mathbb{R}^n$ whose value at x is $\chi(\rho^{-1}|x|)$. For $x \in B$, introduce as shorthand $\mathfrak{h}(t, x) = E_A(t; x, 0)$. Use h to denote $\chi_\rho \mathfrak{h}$. View h as a function with compact support on B with values in $\text{End}(S \otimes_\mathbb{C} U(k))$. By virtue of (2.11), this function obeys

$$\tfrac{d}{dt} h = \partial_\nu \partial_\nu h + cl(F_{A_F}) \chi_{4\rho} h + \chi_\rho V h + \chi_\rho V(1-\chi_\rho)\mathfrak{h}.$$

(2.18)

Of interest here is the trace of $cl(\hat{a}) h$ in $S \otimes_\mathbb{C} \mathbb{C}^k$ at $x = 0$.

To obtain a useful approximation for h, introduce the Euclidean heat kernel,

$$K(t; x, y) = (\tfrac{1}{4\pi t})^{n/2} e^{-|x-y|^2/4t} .$$

(2.19)

If $\psi$ is a map from $(0, 1) \times \mathbb{R}^n$ to $S \otimes_{\mathbb{C}} \mathbb{C}^k$ with compact support along $\mathbb{R}^n$, introduce $\mathcal{K}*\psi$ to denote the map whose value at $(t, x)$ is

$$(\mathcal{K}*\psi)(t,x) = \int_0^t (\int_{\mathbb{R}^n} K(s;x,\cdot)\chi_{4\rho}(cl(F_{A_F}) + \chi_\rho V)\psi(s,\cdot)) \, ds ,$$

(2.20)

Set $k_0$ to denote $K(t; x, 0)$, and set

$$k_1 = k_0 + \int_0^t (\int_{\mathbb{R}^n} K(s;x,\cdot)\chi_\rho V(1-\chi_\rho)\mathfrak{h}) ds .$$

(2.21)

By virtue of (2.16), the endomorphism h obeys

$$h = k_1 + \mathcal{K}*k_1 + (\mathcal{K}*)^2 k_1 + \cdots + (\mathcal{K}*)^{(n-1)/2} k_1 + (\mathcal{K}*)^{(n+1)/2} h .$$

(2.22)

The next lemma takes the first step towards the goal of estimating $tr_S(cl(\hat{a})h)|_{x=0}$.

**Lemma 2.5**: *There exists $\kappa > 0$ which is independent of A, r, and t, and it has the following significance: Define $h_0$ by replacing $k_1$ in (2.22) by $k_0$. Then*

$$|h_1 - h_0|_{x=0} \le e^{-\kappa/t} .$$

*Proof of Lemma 2.5*: The proof exploits four of the various ingredients that enter the definition of $\mathcal{K}*$. Here is the first: The operator V can be written so that each derivative operator that appears in V is multiplied by a term that vanishes at $x = 0$. This is to say that when a given derivative $\partial_v$ appears, it appears as $\sigma(x)\partial_v$, with $\sigma$ a smooth function that obeys $\sigma(0) = 0$. For the terms in V with second derivatives, the derivative of $\sigma$ is uniformly bounded. The terms in V with a single derivative have a version of $\sigma$ with $|d\sigma| \le c_0 r$. Here is a second property: For any $m \ge 0$, there exists a constant $c_m$ that is independent of $x, y \in \mathbb{R}^n$ and $t > 0$ such that

- $|(\nabla^m K)(t; x, y)| \le c_m t^{-m/2} K(2t; x, y)$ ,
- $|x|^m t^{-m/2} K(t; 0, x)| \le c_m K(2t; x, 0)$ .

(2.23)

What follows is a third property: For any $m \ge 0$ and $p > 0$, there exists $c_{m,p}$ that is independent of $x \in \mathbb{R}^n$ and $t > 0$ such that

$$\int_0^t (\int_{\mathbb{R}^n} K(t-s;0,y) \, |y|^m K(ps;y,x) dy) \, ds \leq c_{m,p} t^{1+m/2} K((1+p)t;0,x)$$

(2.24)

The final property is a consequence of Proposition 2.1:

$$\int_0^t (\int_{\mathbb{R}^n} K(s;0,y)(1-\chi_\rho)\mathfrak{h}) \, ds \leq e^{-c/t} .$$

(2.25)

where $c > 0$ is independent of A, r and t.

Granted these four properties, integrate by parts so that the derivative operators in the factors of V that appear in

$$(\mathcal{K}*)^j(k_1 - k_0)$$

(2.26)

act to the left. Now perform the indicated integrals starting from the left-most and working towards the right until only the right-most integral has yet to be done. The inequalities in (2.23)-(2.24) can be used to bound each successive integral in an iterative fashion by a constant multiple of the integral on the left side of (2.14). Note in this regard that the constant in each case is bounded by a factor $c_0 \, rt \, (1 + rt)$ for a suitable constant $c_0$. Once this is done, the final integral has the form of that in (2.25) times a factor that is bounded by $c_0 (rt)^j (1 + rt)^j$.

With Lemma 2.5 in hand, the next lemma supplies an upper bound to the norm at $x = 0$ of the right-most term in (2.22).

**Lemma 2.6**: *There exists* $\kappa \geq 1$ *which is independent of* A, r, *and* t, *and it is such that*

$$|(\mathcal{K}*)^{(n+1)/2} h|_{x=0} \leq \kappa \, r^{(n+1)/2} \, t^{1/2} \, e^{\kappa rt} .$$

*Proof of Lemma 2.6*: The strategy is much like that used in the previous proof. To elaborate, once again evaluate the integrals starting with the left-most. Use (2.23) and (2.24) to successively bound each integral. This scheme eventually produces a bound that reads

$$|(\mathcal{K}*)^{(n+1)/2} h|_{x=0} \leq c_0 \, (tr)^{(n-1)/2} (1+rt)^{(n-1)/2} \int_0^t (\int_{\mathbb{R}^n} K(p(t-s);0,y)(r + r^2 |y|^2) |h(s,y)| dy) \, ds .$$

(2.27)

where $p \geq 1$ is a constant that is independent of t, r and A. To continue, use the bound in (2.4) as an upper bound to $|h(s,y)|$ so as to replace (2.27) with

$$|(\mathcal{K}*)^{(n+1)/2}h|_{x=0} \leq c_0 \, (tr)^{(n-1)/2} \, e^{c_0 rt} \int_0^t (\int_{\mathbb{R}^n} K(p(t-s);0,y)(r + r^2 |y|^2) K(s; y,0) \, dy) \, ds \, .$$

(2.28)

Use (2.23) to bound the integral on the right side of (2.28). Doing so gives the bound stated by the lemma.

According to these last two lemmas,

$$|h|_{x=0} = k_0 + \mathcal{K}*k_0 + (\mathcal{K}*)^2 k_0 + \cdots + (\mathcal{K}*)^{(n-1)/2} k_0 + \mathfrak{p},$$

(2.29)

where $|\mathfrak{p}| \leq c_0 \, r^{(n+1)/2} \, t^{1/2} \, e^{c_0 rt}$. The next step is to evaluate the trace of cl(â) with the various terms $\{(\mathcal{K}*)^j k_0\}_{0 \leq j \leq (n-1)/2}$ that appear in (2.29). This is done just as in Chapter 4 of [BGV]. Here, one uses the fact that $|(\nabla_A)^j F_A| \leq r^{1+j/2}$. The result is

$$tr_{\mathbb{S}}(cl(\hat{a})h)_{x=0} = \pi^{1/2} \, t^{-1/2} \, (\tfrac{1}{2\pi i})^{(n+1)/2} *(\Omega_{\hat{A}} \wedge tr_{\mathbb{C}^k} (\hat{a} \wedge ch(F_{A_F}))) + \mathfrak{r}$$

(2.30)

where $|\mathfrak{r}| \leq c_0 \, t^{-1/2} \, r^{(n-1)/2} \, (rt) \, e^{c_0 rt} \, |\hat{a}|$. The approximation in (2.30) with (2.14) establishe what is claimed by Proposition 2.4.

### 3. The proofs of Theorems 1 and 2

This last section uses the results from the previous section to prove Theorems 1 and 2.

***Proof of Theorem 1:*** To start, suppose that $t \in (0, r^{-1/2})$ and $R \geq 1$ have been chosen for use as described in Section 1b. For each $s \in [0, 1]$, let $\mathfrak{n}_s$ denote the maximal number of linearly independent eigenvectors of $D_{A_s}$ whose eigenvalue has absolute value no greater than $Rt^{-1/2}$. Let $\mathfrak{n}$ denote the maximum from the set $\{\mathfrak{n}_s\}_{s \in [0,1]}$. According to Proposition 2.2, this number $\mathfrak{n}$ obeys $\mathfrak{n} \leq \kappa \, R^n t^{-n/2}$.

For $s \in [0, 1]$, introduce $\wp(s)$ as in (1.13). As explained in Section 1b, the spectral flow $f(A_0, A_1)$ differs from the integral of $\wp$ over the interval $[0, 1]$ by no more than $\mathfrak{n}$. Thus,

$$|f(A_0, A_1) - \int_0^1 \wp(s) \, ds | \leq \kappa \, R^n t^{-n/2} \, .$$

(3.1)

According to Proposition 2.3,

$$\wp(s) = (\tfrac{1}{2\pi i})^{(n+1)/2} \int_M (\Omega_{\hat{A}} \wedge tr_{F \times_{ad} \mathbb{C}^k} (\tfrac{d}{ds} A_{Fs} \wedge ch(F_{A_s}))) + \mathfrak{r}_1$$

(3.2)

where

$$|\mathfrak{r}_1| \leq c_0 (r(A_s)^{(n+1)/2} t + t^{-(n-1)/2} e^{-R^2}) \int_M |\tfrac{d}{ds} A_s| .$$

(3.3)

It is now time to choose values for t and R. To this end, let r denote the maximum from the set $\{r(A_s)\}_{s \in [0,1]}$. Fix $q \in (0, \tfrac{1}{n+1})$, and set $t = r^{-(1+q)}$ and $R = \ln r$. It then follows from (3.1)-(3.3) that $f(A_0, A_1)$ differs from the expression in (1.3) by no more than

$$c_q (1 + (\ln r)^n r^{(n+nq)/2} + \int_0^1 r(A_s)^{(n-1)/2 - q} (\int_M |\tfrac{d}{ds} A_s|) \, ds ))$$

(3.4)

where $c_q$ is a constant that depends on the choice of q. This gives the statement of Theorem 1 with $p = \tfrac{n-1}{2} - q + \varepsilon$ for any sufficiently small but positive $\varepsilon$.

*Proof of Theorem 2*: The proof of Theorem 2 exploits the fact that the difference between the expression in (1.3) and $f(A_0, A_1)$ is gauge invariant. In particular, replace $A_1$ in (1.3) and in $f(A_0, A_1)$ with $A_1' = g \cdot A_1$ where g is a map from M to U(1) with the following properties: First, g has even degree on all generators of $H_1(M; \mathbb{Z})$. Second, $g \cdot A_1 = A_0 + \hat{a}$ where $\hat{a}$ is a coclosed, i-valued 1-form. In addition, the norm of the $L^2$-orthogonal projection of $\hat{a}$ to the space of harmonic 1-forms has norm bounded by $c_0$ where $c_0$ depends only on the Riemannian metric. With $A_1'$ so defined, write $\hat{a} = \hat{a}' + \nu$, where $\nu$ is a harmonic 1-form and where $\hat{a}'$ is $L^2$-orthogonal to the space of harmonic 1-forms. Because $\nu$ is harmonic with uniformly bounded $L^2$ norm, it obeys $|\nu| \leq c_0$.

Now, let G denote the Green's function for the operator $d + d^*$ on the space of forms that are $L^2$-orthogonal to the harmonic forms. Since $d\hat{a}' = F_{A_1} - F_{A_0}$ and $d^*\hat{a}' = 0$, the 1-form $\hat{a}'$ is equal to

$$\hat{a}'|_x = \int_M G_x (F_{A_1} - F_{A_0}) ,$$

(3.4)

Here $G_x(\cdot)$ is smooth where except at x, and near x it obeys $|G_x(y)| \leq c_0 (\tfrac{1}{\text{dist}(x,y)})^{n-1}$. In particular, $|G_x|$ is integrable and so $|\hat{a}'| \leq c_0 r$. Since $|\nu| \leq c_0$, it then follows that $|\hat{a}| \leq c_0 r$ as well.

Granted this estimate, take the path $s \to A_s$ where $A_s = A_0 + s\hat{a}$ and invoke Theorem 1 to obtain Theorem 2.